# ASYMPTOTIC AMPLITUDES AND CAUCHY GAINS: A SMALL-GAIN PRINCIPLE AND AN APPLICATION TO INHIBITORY BIOLOGICAL FEEDBACK

E.D. Sontag[*]
Department of Mathematics
Rutgers University, New Brunswick, NJ 08903
http://www.math.rutgers.edu/~sontag

**Abstract**

The notions of asymptotic amplitude for signals, and Cauchy gain for input/output systems, and an associated small-gain principle, are introduced. These concepts allow the consideration of systems with multiple, and possibly feedback-dependent, steady states. A Lyapunov-like characterization allows the computation of gains for state-space systems, and the formulation of sufficient conditions insuring the lack of oscillations and chaotic behaviors in a wide variety of cascades and feedback loops. An application in biology (MAPK signaling) is worked out in detail.

## 1 Introduction

In this note, we introduce the notions of *asymptotic amplitude* for signals and associated *Cauchy gains* for input/output systems. We present a simple small-gain principle for Cauchy gains, and a Lyapunov-like characterization which allows the estimation of gains for state-space systems.

The concepts and results given here should be of general interest in nonlinear stability and control, especially in those cases in which classical small-gain theorems cannot be applied because the location of closed-loop steady-states depends on the precise gain of the feedback law, or because there are multiple such states.

In developing these ideas, we were originally motivated by the problem of guaranteeing the non-existence of oscillations in certain biological inhibitory feedback loops, and specifically in a mathematical model of *mitogen-activated protein kinase (MAPK) cascades*, which represent a "biological module" or subcircuit which is ubiquitous in eukaryotic cell signal transduction processes. (We are greatly indebted to Stas Shvartsman for bringing to our attention this problem and especially Kholodenko's paper [8], which dealt with the onset of oscillations under high gains.) The general results are illustrated with a numerical computation involving MAPK cascades.

[*]Supported in part by US Air Force Grant F49620-01-1-0063

## 1.1 Cauchy Gains

For any metric space $M$, we write the distance $d_M(a,b)$ between any two elements $a, b \in M$, in the suggestive form "$|a - b|$" even when $M$ has no linear structure (so the "$-$" sign has no meaning, of course), and define the *asymptotic amplitude* of a function $\omega : \mathbb{R}_{\geq 0} \to M$, where $\mathbb{R}_{\geq 0} = [0, +\infty)$, as follows:

$$\|\omega\|_{\mathrm{aa}} := \limsup_{s,t \to \infty} |\omega(t) - \omega(s)| = \lim_{T \to \infty} \left( \sup_{t,s \geq T} |\omega(t) - \omega(s)| \right) \in [0, \infty].$$

Observe that the condition "$\|\omega\|_{\mathrm{aa}} = 0$" amounts to the Cauchy property for $\omega$: for every $\varepsilon > 0$ there is some $T > 0$ such that $|\omega(t) - \omega(s)| < \varepsilon$ for all $t, s \geq T$. Thus, when $M$ is a complete metric space (for instance, if, as in all our examples, $M \subseteq \mathbb{R}^m$ is any closed subset of a Euclidean space):

$$\|\omega\|_{\mathrm{aa}} = 0 \iff \exists \lim_{t \to \infty} \omega(t).$$

If $\|\omega\|_{\mathrm{aa}} = 0$, we denote $\omega^{\infty} := \lim_{t \to \infty} \omega(t)$.

Let $\mathcal{U}$ and $\mathcal{Y}$ be two complete metric spaces. We define a *behavior with input-value space $\mathcal{U}$ and output-value space $\mathcal{Y}$* as a relation $\mathcal{R}$ between time-functions with values in $\mathcal{U}$ and $\mathcal{Y}$ respectively:

$$\mathcal{R} \subseteq [\mathbb{R}_{\geq 0} \to \mathcal{U}] \times [\mathbb{R}_{\geq 0} \to \mathcal{Y}]$$

where $[\mathbb{R}_{\geq 0} \to M]$ is the set of functions $\mathbb{R}_{\geq 0} \to M$. We call any element $(\omega, \eta) \in \mathcal{R}$ an *input/output pair*, and say that $\omega$ is an input signal and $\eta$ is an output signal of $\mathcal{R}$.

Typical examples of behaviors, to be discussed in detail later, are those obtained by starting with a system of differential equations with inputs ("forcing functions" or "controls") $\omega$, and viewing the solutions obtained by solving the system with different initial states, or some components of these solutions, as the outputs $\eta$. The formalism that we use, based on abstract relations for the formulation of small-gain results, dates back to Zames' original paper [21], and the term "behavior" is borrowed from Willems' work [20].

We use standard terminology for comparison functions: $\mathcal{K}_{\infty}$ is the class of continuous, strictly increasing, and unbounded functions $\gamma : \mathbb{R}_{\geq 0} \to \mathbb{R}_{\geq 0}$ with $\gamma(0) = 0$.

**Definition 1.1** A behavior $\mathcal{R}$ has *Cauchy gain* $\gamma \in \mathcal{K}_{\infty}$ if

$$\|\eta\|_{\mathrm{aa}} \leq \gamma(\|\omega\|_{\mathrm{aa}})$$

for all $(\omega, \eta) \in \mathcal{R}$. □

The existence of a Cauchy gain for $\mathcal{R}$ implies, in particular, the following *converging input converging output* property for $\mathcal{R}$: if $\omega(t) \to \bar{u}$ as $t \to \infty$, for some $\bar{u} \in \mathcal{U}$ (that is, if $\|\omega\|_{\mathrm{aa}} = 0$), and if $(\omega, \eta) \in \mathcal{R}$, then also $\eta(t) \to \bar{y}$ as $t \to \infty$, for some $\bar{y} \in \mathcal{Y}$.

The interconnection that results when the output of a system $\mathcal{R}$ is fed back to its input under the action of the system (feedback law) $\mathcal{S}$ is pictorially represented in Figure 1. The behavioral terminology gives an easy way to define formally the meaning of this interconnection: if $\mathcal{R}$ and $\mathcal{S}$ are behaviors, then the signals that appear when the loop is closed are precisely those pairs $(\omega, \eta)$ such that $(\omega, \eta) \in \mathcal{R}$ and $(\eta, \omega) \in \mathcal{S}$. Put another way, the feedback connection is simply the behavior $\mathcal{R} \bigcap \mathcal{S}^{-1}$, where, for any behavior $\mathcal{S}$ with input-value space $\mathcal{Y}$ and output-value space $\mathcal{U}$, we denote by $\mathcal{S}^{-1}$ the inverse behavior, with input-value space $\mathcal{U}$ and output-value space $\mathcal{Y}$, consisting of all pairs $(\omega, \eta)$ such that $(\eta, \omega) \in \mathcal{S}$.

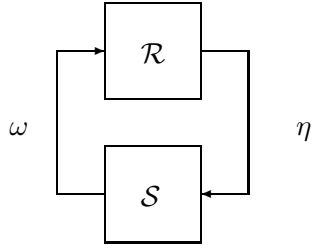

Figure 1: Feedback Interconnection $\mathcal{R} \bigcap \mathcal{S}^{-1}$

With this formalism, the basic "small gain principle" is trivial to establish. It states that the interconnection of two systems having Cauchy gains whose composition is a contraction, has the property that the external signals $\omega$ and $\eta$ must always converge to some value as $t \to \infty$, at least if they are known to have finite asymptotic amplitude:

**Lemma 1.2** *(Small gain lemma for asymptotic amplitude.)* Suppose that $\mathcal{R}$ and $\mathcal{S}$ are two behaviors with Cauchy gains $\gamma_1$ and $\gamma_2$ respectively, and that the following condition holds:

$$\gamma_1(\gamma_2(r)) < r \quad \forall\, r > 0\,. \tag{1}$$

Then, for all $(\omega, \eta) \in \mathcal{R} \bigcap \mathcal{S}^{-1}$ for which $\|\omega\|_{\mathrm{aa}} < \infty$, $\|\omega\|_{\mathrm{aa}} = \|\eta\|_{\mathrm{aa}} = 0$.

*Proof.* Since $(\omega, \eta) \in \mathcal{R}$, $\|\eta\|_{\mathrm{aa}} \leq \gamma_1(\|\omega\|_{\mathrm{aa}})$; and since also $(\eta, \omega) \in \mathcal{S}$, $\|\omega\|_{\mathrm{aa}} \leq \gamma_2(\|\eta\|_{\mathrm{aa}})$. If $\|\eta\|_{\mathrm{aa}} \neq 0$, then $\|\eta\|_{\mathrm{aa}} \leq \gamma_1(\gamma_2(\|\eta\|_{\mathrm{aa}})) < \|\eta\|_{\mathrm{aa}}$, a contradiction. Finally, $\|\omega\|_{\mathrm{aa}} \leq \gamma_2(\|\eta\|_{\mathrm{aa}}) = \gamma_2(0) = 0$ gives that also $\|\omega\|_{\mathrm{aa}} = 0$. ∎

**Remark 1.3** Note that the condition "$\|\omega\|_{\mathrm{aa}} < \infty$" is equivalent to ultimate boundedness, i.e. there are a bounded set $C \subseteq \mathcal{U}$ and some $T \geq 0$ such that $\omega(t) \in C$ for all $t \geq T$. (Writing $|u| := |u - 0|$ for some fixed element $0 \in \mathcal{U}$: if there are some $c, T > 0$ so that $|\omega(t)| \leq c$ for all $t \geq T$ then $\|\omega\|_{\mathrm{aa}} \leq 2c$; conversely, if $\sup_{t,s \geq T} |\omega(t) - \omega(s)| \leq c$ for some $T$ then $|w(t)| \leq c + |\omega(T)|$ for all $t \geq T$.) In applications to feedback loops involving differential equations, all signals are continuous, and for them, ultimate boundedness is equivalent to just boundedness. □

The limiting values of the signals $\omega$ and $\eta$, whose existence is asserted by Lemma 1.2, need not be unique; for instance bistable systems give rise to nonuniqueness. In order to present a condition which guarantees uniqueness, we introduce a new concept.

**Definition 1.4** A behavior $\mathcal{R}$ has *incremental limit gain* $\kappa \in \mathcal{K}_\infty$ if the following property holds:
$$\limsup_{t \to \infty} |\eta_1(t) - \eta_2(t)| \leq \kappa(|\omega_1^\infty - \omega_2^\infty|)$$
whenever $(\omega_i, \eta_i) \in \mathcal{R}$ are any two pairs with the properties $\|\omega_1\|_{\mathrm{aa}} = \|\omega_2\|_{\mathrm{aa}} = 0$. □

In words, this definition says that, if we are given two input/output pairs for which the inputs converge, and if the limits of the two inputs are close to each other, then the outputs become asymptotically close to each other. If $\mathcal{R}$ has an incremental limit gain $\kappa$, and if in addition $\mathcal{R}$ also admits a Cauchy gain, then both $\eta_1^\infty$ and $\eta_2^\infty$ exist whenever $\|\omega_1\|_{\mathrm{aa}} = \|\omega_2\|_{\mathrm{aa}} = 0$ (converging-input converging-output), and thus the "limsup" in Definition 1.4 is a limit, and the estimate becomes:

$$|\eta_1^\infty - \eta_2^\infty| \leq \kappa(|\omega_1^\infty - \omega_2^\infty|). \tag{2}$$

With this concept, we have another obvious observation:

**Lemma 1.5** *(Small gain lemma for asymptotic amplitude, with uniqueness.)* Suppose that $\mathcal{R}$ and $\mathcal{S}$ are two behaviors with Cauchy gains $\gamma_1$ and $\gamma_2$ respectively, and incremental limit gains $\kappa_1$ and $\kappa_2$ respectively, and that the following condition holds:

$$\kappa_1(\kappa_2(r)) < r \quad \forall\, r > 0 \tag{3}$$

in addition to (1). Then, there exist two elements $\bar{u} \in \mathcal{U}$ and $\bar{y} \in \mathcal{Y}$ such that, for every input/output pair $(\omega, \eta) \in \mathcal{R} \bigcap \mathcal{S}^{-1}$ for which $\|\omega\|_{\mathrm{aa}} < \infty$, $\omega^\infty = \bar{u}$ and $\eta^\infty = \bar{y}$.

*Proof.* If $\mathcal{R} \bigcap \mathcal{S}^{-1} = \emptyset$, there is nothing to prove. Otherwise, pick an arbitrary $(\omega_1, \eta_1) \in \mathcal{R} \bigcap \mathcal{S}^{-1}$ for which $\|\omega_1\|_{\mathrm{aa}} < \infty$. From Lemma 1.2, there exist $\bar{u} := \omega_1^\infty$ and $\bar{y} := \eta_1^\infty$. Pick now any other $(\omega_2, \eta_2) \in \mathcal{R} \bigcap \mathcal{S}^{-1}$ for which $\|\omega_2\|_{\mathrm{aa}} < \infty$; again by the Lemma, $\omega_2^\infty$ and $\eta_2^\infty$ exist. By the incremental limit gain property, in the form (2), both $|\bar{y} - \eta_2^\infty| \leq \kappa_1(|\bar{u} - \omega_2^\infty|)$ and $|\bar{u} - \omega_2^\infty| \leq \kappa_2(|\bar{y} - \eta_2^\infty|)$. From

$$|\bar{y} - \eta_2^\infty| \leq \kappa_1(\kappa_2(|\bar{y} - \eta_2^\infty|))$$

we conclude that $\eta_2^\infty = \bar{y}$, and so also $\omega_2^\infty = \bar{u}$. ∎

Once the appropriate definitions have been given, the two Lemmas are quite obvious. The harder step is, often, to verify when the Lemmas apply. In order to carry out such an application, one needs to find sufficient and easy to check conditions which guarantee the existence of Cauchy and incremental limit gains, for the systems whose feedback interconnection is being studied.

We will mainly study behaviors $\mathcal{R}$ which can be built up from cascades of simpler behaviors $\mathcal{R}_i$, each of which is either defined by a system of differential equations, by a pure delay, or by a memoryless nonlinearity. The composition $\mathcal{R}$ will represent the input/output pairs of a large set of delay-differential equations. The Cauchy and incremental limit gains of the behaviors $\mathcal{R}_i$ can be composed, so as to provide the gains of the complete system $\mathcal{R}$. Section 2 describes these general ideas. Section 3 shows how to estimate gains based on contractions of omega-limit sets of signals, and these types of estimates are used in order to justify the Lyapunov-like methods described in Section 4 for state-space systems. Finally, Section 5 specialized the results to a class of inhibitory feedback loops, and in particular for the motivating MAPK example.

This work is related to other work on "nonlinear gain" small-gain theorems such as [7, 9, 5], which in turn was motivated by classical small-gain theorems as in [1, 11, 12, 21]. Future developments will include generalizations to estimates which quantify overshoot, in the ISS (cf. [15]) sense.

## 2 Simple Behaviors and Cascades

The *delay-$\tau$* operator $\mathcal{D}_\tau$ on $\mathcal{U}$, where $\tau \geq 0$, is the behavior, with $\mathcal{Y} = \mathcal{U}$, defined by: $(\omega, \eta) \in \mathcal{D}_\tau$ if and only if $\eta(t) = \omega(t - \tau)$ for all $t \geq \tau$. (The value of the output for $t < \tau$ is arbitrary; in an abstract dynamical systems sense, it forms part of the specification of initial conditions.) It is clear that $\mathcal{D}_\tau$ has Cauchy gain $I$ and incremental limit gain $I$, where $I$ is the identity function, $I(r) = r$.

Given any map $\psi : \mathcal{U} \to \mathcal{Y}$, the *memoryless behavior associated to $\psi$*, which we denote by $\mathcal{M}_\psi$, is the behavior consisting of all pairs of functions $(\omega, \eta)$ such that $\eta(t) = \psi(\omega(t))$ for all $t$. Suppose that $\psi$ is a Lipschitz map: for some $\lambda \geq 0$, $|\psi(u_1) - \psi(u_2)| \leq \lambda |u_1 - u_2|$ for all $u_1, u_2 \in \mathcal{U}$. Then $\mathcal{M}_\psi$ has Cauchy gain $\lambda I$ and incremental limit gain $\lambda I$, where $\lambda I(r) = \lambda r$.

Suppose that $\mathcal{R} \subseteq [\mathbb{R}_{\geq 0} \to \mathcal{U}] \times [\mathbb{R}_{\geq 0} \to \mathcal{Y}]$ and $\mathcal{S} \subseteq [\mathbb{R}_{\geq 0} \to \mathcal{Y}] \times [\mathbb{R}_{\geq 0} \to \mathcal{Z}]$ are two behaviors, with Cauchy gains $\gamma_1$ and $\gamma_2$ respectively, and consider the cascade combination shown pictorially in Figure 2 and defined formally as:

$$\mathcal{S} \circ \mathcal{R} := \{(\omega, \zeta) \mid (\exists \eta \in [\mathbb{R}_{\geq 0} \to \mathcal{Y}]) \text{ s.t. } (\omega, \eta) \in \mathcal{R} \ \& \ (\eta, \zeta) \in \mathcal{S}\} .$$

Then, clearly, $\mathcal{S} \circ \mathcal{R}$ has Cauchy gain $\gamma_2 \circ \gamma_1$. Suppose now that also $\mathcal{R}$ and $\mathcal{S}$ have incremental limit gains $\kappa_1$ and $\kappa_2$ respectively. Let $(\omega_i, \eta_i) \in \mathcal{R}$ and $(\eta_i, \zeta_i) \in \mathcal{S}$, $\|\omega_i\|_{\text{aa}} = 0$, for $i = 1, 2$. We have that $\eta_1^\infty$ and $\eta_2^\infty$ exist, and (2) holds with $\kappa = \kappa_1$. Similarly, since $\mathcal{S}$ has a Cauchy gain, $\zeta_1^\infty$ and $\zeta_2^\infty$ exist, and $|\zeta_1^\infty - \zeta_2^\infty| \leq \kappa_2(|\eta_1^\infty - \eta_2^\infty|)$. Therefore

$$|\zeta_1^\infty - \zeta_2^\infty| \leq \kappa_2(\kappa_1(|\omega_1^\infty - \omega_2^\infty|))$$

and hence $\mathcal{S} \circ \mathcal{R}$ has incremental limit gain $\kappa_2 \circ \kappa_1$.

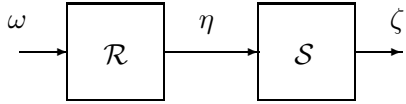

Figure 2: Cascade $\mathcal{S} \circ \mathcal{R}$

### 2.1 Tighter Estimates: Relative Gains

Tighter estimates of gains for the cascade $\mathcal{S} \circ \mathcal{R}$ can make use of the following observation. Suppose that the possible output signals of $\mathcal{R}$ all tend, as $t \to \infty$, to values in a restricted subset $Y$ of $\mathcal{Y}$. Then the relevant gains $\gamma_2, \kappa_2$ should be the gains of $\mathcal{S}$ when restricted to those signals of the form $(\eta, \zeta) \in \mathcal{S}$ such that $\eta \in [\mathbb{R}_{\geq 0} \to Y]$. These gains may well be smaller than the original ones, so that smaller overall gains result for the cascade. Let us make this precise.

For any subset $U_0 \subseteq \mathcal{U}$, we write "$\omega \to U_0$" if $\omega(t)$ converges to $U_0$ as $t \to \infty$, that is, for every $\varepsilon > 0$ there is some $T \geq 0$ such that

$$\omega(t) \in B_\varepsilon(U_0) = \{u \in \mathcal{U} \mid (\exists u' \in U_0) |u - u'| \leq \varepsilon\}$$

for every $t \geq T$.

Let $U_0 \subseteq \mathcal{U}$ and let $\mathcal{R} \subseteq [\mathbb{R}_{\geq 0} \to \mathcal{U}] \times [\mathbb{R}_{\geq 0} \to \mathcal{Y}]$. We will say that $\mathcal{R}$ has a *Cauchy gain $\gamma$ on $U_0$* if $\|\eta\|_{\text{aa}} \leq \gamma(\|\omega\|_{\text{aa}})$ holds for each input/output pair $(\omega, \eta) \in \mathcal{R}$ for which $\omega \to U_0$.

Similarly, we say that $\mathcal{R}$ has *incremental limit gain* $\kappa$ on $U_0$ if $\limsup_{t\to\infty} |\eta_1(t) - \eta_2(t)| \leq \kappa(|\omega_1^\infty - \omega_2^\infty|)$ holds whenever $(\omega_i, \eta_i) \in \mathcal{R}$ are any two pairs such that $\omega_1^\infty$ and $\omega_2^\infty$ both exist and belong to $U_0$. In the special case $U_0 = \mathcal{U}$, one recovers the definitions of Cauchy and incremental limit gains.

Suppose now that there are two sets $U_0 \subseteq \mathcal{U}$ and $Y_0 \subseteq \mathcal{Y}$ such that:

- $\mathcal{R}$ has Cauchy gain $\gamma_1$ on $U_0$.
- $\mathcal{S}$ has Cauchy gain $\gamma_2$ on $Y_0$.
- Whenever $(\omega, \eta) \in \mathcal{R}$ is so that $\omega \to U_0$, necessarily $\eta \to Y_0$.

Then, clearly, $\mathcal{S} \circ \mathcal{R}$ has Cauchy gain $\gamma_2 \circ \gamma_1$ on $U_0$. An analogous conclusion holds for incremental limit gain on $U_0$.

## 3 A Sufficient Condition

Recall that, for any metric space $\mathcal{U}$ and function $\omega : \mathbb{R}_{\geq 0} \to \mathcal{U}$, the omega-limit set $\Omega = \Omega^+[\omega]$ is the set consisting of those points $u \in \mathcal{U}$ for which there exists a convergent sequence $\omega(t_i) \to u$, for some sequence $\{t_i\} \subseteq \mathbb{R}_{\geq 0}$ such that $t_i \to \infty$ as $i \to \infty$. The following properties are elementary: (1) the set $\Omega$ is closed; (2) if $\omega \to U$ and $U$ is closed, then $\Omega \subseteq U$, and so $\Omega$ is compact if $U$ is; (3) provided $\omega$ is precompact, that is to say, if there is some compact subset $U \subseteq \mathcal{U}$ such that $\omega(t) \in U$ for all $t \geq 0$, then $\Omega$ is compact, and $\omega \to \Omega$ (proof of this last statement: if there is some $\varepsilon > 0$ and some sequence $t_i \to \infty$ such that $\omega(t_i) \in U \setminus B_\varepsilon(\Omega)$ for all $i$, then one can pick a subsequence of $\{t_i\}$ such that $\omega(t_{i_j}) \to u$ for some $u$, and thus $u \in U \setminus B_{\varepsilon/2}(\Omega)$, a contradiction since $u \in \Omega$ by definition of $\Omega$); (4) if $\omega$ is precompact and is a continuous function of $t$, then $\Omega$ is connected (proof: otherwise, there are two nonempty compact subsets with $\Omega = \Omega_1 \bigcup \Omega_2$ and, for some $\varepsilon > 0$, $\mathrm{dist}(\Omega_1, \Omega_2) > 3\varepsilon$, and by (3) some $T > 0$ such that $\omega(t) \in B = B_\varepsilon(\Omega)$ for all $t \geq T$; from $\Omega_1 \subseteq \Omega$ we know that there is some $t_1 > T$ such that $\omega(t_1) \in B_1 = B_\varepsilon(\Omega_1)$ and from $\Omega_2 \subseteq \Omega$ that there is a $t_2 > t_1$ such that $\omega(t_2) \in B_2 = B_\varepsilon(\Omega_2)$, so, noticing that $B = B_1 \bigcup B_2$ and that $B_1 \bigcap B_2 = \emptyset$ by choice of $\varepsilon$, and writing $\mathcal{I} := \{\omega(t), t \in [t_1, t_2]\}$, we have that this connected set can be written as a disjoint union $\mathcal{I} = (\mathcal{I} \bigcap B_1) \bigcup (\mathcal{I} \bigcap B_2)$ of nonempty closed sets, a contradiction).

In general, we denote by $|U|$ the *diameter* $\sup\{|u - v| \mid u, v \in U\}$ of a closed subset $U$ of a metric space $\mathcal{U}$. For each $\omega : \mathbb{R}_{\geq 0} \to \mathcal{U}$, it holds that $|\Omega^+[\omega]| \leq \|\omega\|_{\mathrm{aa}}$, and, if $\omega$ is precompact,

$$\|\omega\|_{\mathrm{aa}} = |\Omega^+[\omega]|. \tag{4}$$

Indeed, pick any $\varepsilon > 0$ and two elements $u, v \in \Omega$ such that $|u - v| \geq |\Omega| - \varepsilon$; then there are two sequences $\omega(t_i) \to u$ and $\omega(s_i) \to v$, so $\|\omega\|_{\mathrm{aa}} = \limsup_{s,t\to\infty} |\omega(t) - \omega(s)| \geq |u - v| \geq |\Omega| - \varepsilon$. As this is true for every $\varepsilon > 0$, we have $\|\omega\|_{\mathrm{aa}} \geq |\Omega|$. Conversely, if $|\omega(t_i) - \omega(s_i)| \geq \|\omega\|_{\mathrm{aa}} - \varepsilon$ for some two sequences $t_i \to \infty$ and $s_i \to \infty$, we may extract first a subsequence of $\{t_i\}$ such that $\omega(t_{i_j})$ is convergent (precompactness is used here), and then a subsequence of $\{s_{i_j}\}$, so that, without loss of generality we may suppose that $\omega(t_i) \to u$ and $\omega(s_i) \to v$ for some $u, v \in \Omega$, and thus $|\Omega| \geq |u - v| \geq \|\omega\|_{\mathrm{aa}} - \varepsilon$, so letting $\varepsilon \to 0$ gives the other inequality.

We say that a mapping $\Gamma$ assigning subsets of one set to subsets of another is monotonic if $U_1 \subseteq U_2 \Rightarrow \Gamma(U_1) \subseteq \Gamma(U_2)$.

**Lemma 3.1** Suppose given a behavior $\mathcal{R}$, a compact subset $U_0 \subseteq \mathcal{U}$, a function $\gamma \in \mathcal{K}_\infty$, and a mapping $\Gamma$ from compact subsets of $U_0$ to subsets of $\mathcal{Y}$, such that the following properties hold:

(a) For each $(\omega, \eta) \in \mathcal{R}$ for which $\Omega^+[\omega] \subseteq U_0$, the output $\eta$ is precompact.

(b) For each compact subset $U \subseteq U_0$, and each $(\omega, \eta) \in \mathcal{R}$ for which $\Omega^+[\omega] \subseteq U$, it holds that $\Omega^+[\eta] \subseteq \Gamma(U)$.

(c) For each compact subset $U \subseteq U_0$, it holds that $|\Gamma(U)| \leq \gamma(|U|)$.

Then $\mathcal{R}$ has Cauchy gain $\gamma$ on $U_0$ and incremental limit gain $\gamma$ on $U_0$. Moreover, for each compact subset $U \subseteq U_0$, and each $(\omega, \eta) \in \mathcal{R}$ for which $\omega \to U$, $\eta \to \Gamma(\Omega^+[\omega])$. If $\Gamma$ is monotonic, then also $\eta \to \Gamma(U)$.

*Proof.* Pick any $(\omega, \eta) \in \mathcal{R}$ and any compact $U \subseteq U_0$, and suppose that $\omega \to U$. By (1) and (2) in the previous discussion, the set $\Omega^+[\omega]$ is a compact subset of $U$. By (a), $\eta$ is precompact. Therefore $\|\eta\|_{\text{aa}} = |\Omega^+[\eta]|$, and also $\eta \to \Omega^+[\eta]$. By (b), applied to $\Omega^+[\omega]$ itself, we know that $\Omega^+[\eta] \subseteq \Gamma(\Omega^+[\omega])$, which gives the conclusion $\eta \to \Gamma(\Omega^+[\omega])$. If $\Gamma$ is monotonic, then $\Omega^+[\omega] \subseteq U$ implies that $\Gamma(\Omega^+[\omega]) \subseteq \Gamma(U)$, so $\eta \to \Gamma(U)$. In addition, $|\Omega^+[\eta]| \leq |\Gamma(\Omega^+[\omega])|$ together with (c) give the following inequality:

$$\|\eta\|_{\text{aa}} = |\Omega^+[\eta]| \leq |\Gamma(\Omega^+[\omega])| \leq \gamma(|\Omega^+[\omega]|) \leq \gamma(\|\omega\|_{\text{aa}}). \tag{5}$$

When applied in the special case $U = U_0$, this establishes the Cauchy gain conclusion.

Suppose now that $(\omega_i, \eta_i) \in \mathcal{R}$ are any two pairs such that $\omega_1^\infty$ and $\omega_2^\infty$ both exist and belong to $U_0$. In particular, $\omega_1 \to U_0$ and $\omega_2 \to U_0$. So both $\eta_1^\infty$ and $\eta_2^\infty$ exist, by the Cauchy gain conclusion. Note that $\Omega^+[\omega_i] = \{\omega_i^\infty\}$ and $\Omega^+[\eta_i] = \{\eta_i^\infty\}$ for $i = 1, 2$. We introduce the two-element set $U = \{\omega_1^\infty, \omega_2^\infty\} \subseteq U_0$; note that $|U| = |\omega_1^\infty - \omega_2^\infty|$. From $\Omega^+[\omega_i] \subseteq U$ and (b), we have that $\Omega^+[\eta_i] \subseteq \Gamma(U)$, that is, $\eta_i^\infty \in \Gamma(U)$, for $i = 1, 2$. Therefore

$$|\eta_1^\infty - \eta_2^\infty| \leq |\Gamma(U)| \leq \gamma(|U|) = \gamma(|\omega_1^\infty - \omega_2^\infty|),$$

which proves the incremental limit property. ∎

**Remark 3.2** The Cauchy gain (but not the incremental limit gain) results in Lemma 3.1 can be tightened provided that one knows that $\eta$ is continuous, as is the case when considering behaviors defined by differential equation systems, and even further provided that $\omega$ is continuous, as is the case when dealing with feedback configurations $\mathcal{R} \cap \mathcal{S}^{-1}$, in the latter case assuming that $\mathcal{U}$ is locally compact (which is automatically satisfied in all finite dimensional applications of the results). To be precise, let us denote by $|U|_c$ the maximal diameter of the connected components of a set $U$, that is, the supremum of the quantities $|u - v|$ taken over all pairs $u, v$ which lie in any given connected component of $U$. In general, for any continuous and precompact $\omega$, and any subset $U \subseteq \mathcal{U}$ such that $\Omega^+[\omega] \subseteq U$, it holds that $\|\omega\|_{\text{aa}} \leq |U|_c$, because $\Omega^+[\omega]$ is connected and thus lies entirely in a single connected component of $U$. Then:

- if $\eta$ is continuous for every $(\omega, \eta) \in \mathcal{R}$ then assumption (c) can be weakened to "$|\Gamma(U)|_c \leq \gamma(|U|)$ for every compact subset $U \subseteq U_0$",

- if in addition $\omega$ is continuous for every $(\omega, \eta) \in \mathcal{R}$ and $\mathcal{U}$ is locally compact then assumption (c) can be weakened to "$|\Gamma(U)|_c \leq \gamma(|U|)$ for every compact and connected subset $U \subseteq U_0$".

Indeed, in the proof of Lemma 3.1 observe that if $\eta$ is precompact and continuous then $\Omega^+[\eta]$ is connected, so $\Omega^+[\eta] \subseteq \Gamma(\Omega^+[\omega])$ implies that $\|\eta\|_{\mathrm{aa}} \leq |\Gamma(\Omega^+[\omega])|_c$, and this is all that is needed in (5). If also $\mathcal{U}$ is locally compact and $\omega$ is continuous, then $\omega \to U$ (compact) implies that $\omega$ is precompact, which together with continuity says that $\Omega^+[\omega]$ is connected, and therefore in order to prove (5), (c) only needs to be applied with the connected set $U = \Omega^+[\omega]$. □

## 4  Systems of Differential Equations

A particular class of behaviors, in fact the main objects of interest in this note, are obtained as follows. We consider systems of differential equations with inputs and outputs:

$$\dot{x} = f(x, u), \quad y = h(x) \tag{6}$$

for which states $x(t)$ evolve in a subset $\mathcal{X}$ of a Euclidean space $\mathbb{R}^n$, inputs take values $u(t)$ in a complete metric space $\mathcal{U}$ and outputs take values $y(t)$ in a complete metric space $\mathcal{Y}$. (Typically in applications, $\mathcal{U}$ and $\mathcal{Y}$ are any two closed subsets of Euclidean spaces.) Technically, we assume that $f : \mathcal{X}_0 \times \mathcal{U} \to \mathbb{R}^n$ is defined on an open subset $\mathcal{X}_0 \subseteq \mathbb{R}^n$ which contains $\mathcal{X}$, is continuous, and is locally Lipschitz in $x$ uniformly on compact subsets of $\mathcal{X}_0 \times \mathcal{U}$; the map $h : \mathcal{X} \to \mathcal{Y}$ is assumed to be continuous. Furthermore, $\mathcal{X}$ is an invariant and forward complete subset, in the sense that, for each Lebesgue-measurable precompact input $\omega : \mathbb{R}_{\geq 0} \to \mathcal{U}$, and each initial state $x_0 \in \mathcal{X}$, the unique solution $\xi(t) = \varphi(t, x_0, \omega)$ of the initial value problem $\dot{\xi}(t) = f(\xi(t), \omega(t))$, $\xi(0) = x_0$, is defined and satisfies $\xi(t) \in \mathcal{X}$ for all $t \geq 0$. (The function $\xi$ is Lipschitz, and hence differentiable almost everywhere; if $\omega$ is continuous, then $\xi$ is continuously differentiable.) To any given system (6) one associates a behavior $\mathcal{R}$, with input-value space $\mathcal{U}$ and output-value space $\mathcal{Y}$, defined by: $(\omega, \eta) \in \mathcal{R}$ if and only if $\omega$ is precompact and Lebesgue-measurable, and there exists some $x_0 \in \mathcal{X}$ such that $\eta(t) = h(\varphi(t, x_0, \omega))$ for all $t \in \mathbb{R}_{\geq 0}$. We call $\mathcal{R}$ *the behavior of* (6).

**Remark 4.1** A minor technicality concerns the fact that Lebesgue-measurable functions are, strictly speaking, not functions but equivalence classes of functions, so one should interpret the "limsup" in the definition of asymptotic amplitude in an "almost everywhere" manner; similarly, we interpret "precompact" as meaning that there is some $\omega$ in the given equivalence class whose values all remain in a compact. From now on, we leave this technicality implicit; in applications to stability of feedback loops involving systems of differential equations, all the functions considered are continuous –even differentiable– so the issue does not even arise. □

We will obtain sufficient conditions for the existence of the two types of gains, expressed in terms of Lyapunov-type functions.

Given a subset $U \subseteq \mathcal{U}$, we will say that a function

$$V : \mathcal{X} \to \mathbb{R}_{\geq 0}$$

is a *U-decrease* function provided that the following properties hold:

- $V$ is proper, that is, $\{x \in \mathcal{X} \mid a \leq V(x) \leq b\}$ is a compact subset of $\mathcal{X}$, for each $a \leq b$ in $\mathbb{R}_{\geq 0}$;

- $V$ is continuous;

- for each $x \in \mathcal{X}$ which does not belong to $Z_V := \{x \mid V(x) = 0\}$, it holds that $V$ is continuously differentiable in a neighborhood of $x$ and

$$\nabla V(x) \cdot f(x, u) < 0 \qquad (7)$$

for all $u \in U$.

(We understand continuous differentiability in the following sense: there is a neighborhood of $x$ in $\mathcal{X}_0$ such that $V$ extends to this neighborhood as a $\mathcal{C}^1$ function.)

**Lemma 4.2** Suppose that $V$ is a $U$-decrease function, for some compact subset $U \subseteq \mathcal{U}$. Pick any Lebesgue-measurable precompact $\omega : \mathbb{R}_{\geq 0} \to \mathcal{U}$ and any solution $\xi$ of the system $\dot{\xi} = f(\xi, \omega)$. Suppose that either:

1. there is some $T_0 \geq 0$ such that $\omega(t) \in U$ for all $t \geq T_0$, or

2. $\omega \to U$ and $\xi$ is precompact.

Then $\xi \to Z_V$.

*Proof.* Given any $\omega$ and $\xi$, we will first let $a > 0$ be arbitrary and prove that the set $V_a := \{x \mid V(x) \leq a\}$ has the property that, for some $T^* \geq 0$,

$$\xi(t) \in V_a \quad \forall t \geq T^*. \qquad (8)$$

If the assumption is that $\omega \to U$ as $t \to \infty$ and that $\xi$ is precompact, that is to say, there is some compact subset $C_0$ of $\mathcal{X}$ such that $\xi(t) \in C_0$ for all $t \geq 0$, then we introduce $b := \max\{V(x), x \in C_0\}$ and the set $C := V^{-1}([a,b])$. Note that $\xi(t) \in C$ whenever $V(\xi(t)) \geq a$, by the choice of $b$. Since, by properness of $V$, $C$ is a compact subset of $\mathcal{X} \setminus Z_V$, by Property (7) there is some $\alpha > 0$ and some neighborhood $\widetilde{U}$ of $U$ in $\mathcal{U}$ so that $\nabla V(x) \cdot f(x, u) \leq -\alpha$ for all $x \in C$ and all $u \in \widetilde{U}$. Since $\omega \to U$, there must be some $T_0$ such that $\omega(t) \in \widetilde{U}$ for all $t \geq T_0$. Thus,

$$\nabla V(x) \cdot f(x, \omega(t)) \leq -\alpha < 0 \qquad \forall x \in C \quad \forall t \geq T_0. \qquad (9)$$

If, instead, the assumption is that there is some $T_0 \geq 0$ such that $\omega(t) \in U$ for all $t \geq T_0$, we pick $T_0$ to be any such number, and let $C := \{x \mid V(x) \geq a\}$. So once more we have that $\nabla V(x) \cdot f(x, \omega(t)) < 0$ for all $x \in C$ and all $t \geq T_0$.

*Claim:* Let $T_1 := \inf\{t \geq T_0 \mid V(\xi(t)) \leq a\}$ (if $V(\xi(t)) > a$ for all $t \geq T_0$, we define $T_1 = +\infty$). Then $V(\xi(t)) \leq V(\xi(T_0))$ for all $t \in [T_0, T_1)$ and $V(\xi(t)) \leq a$ for all $t \geq T_1$.

*Proof of the claim:* Suppose that $\xi(t) \in C$ for all $t$ in some interval $(\tau_1, \tau_2)$ with $\tau_1 \geq T_0$. Then $dV(\xi(t))/dt = \nabla V(\xi(t)) \cdot f(\xi(t), \omega(t)) < 0$ for almost every $t \in (\tau_1, \tau_2)$. Therefore, $V(\xi(t))$ is decreasing on this interval, and we have that $V(\xi(t)) \leq V(\xi(\tau_1))$ for all $t \in (\tau_1, \tau_2)$. In particular, for each $t \in (T_0, T_1)$, by minimality of $T_1$ we know that $V(\xi(t)) > a$ and so $\xi(t) \in C$. This proves the first part of the claim: $V(\xi(t)) \leq V(\xi(T_0))$ for all $t \in [T_0, T_1)$. If $T_1 = \infty$, there is nothing more to prove. So assume that $T_1 < \infty$ and there exists some $S > T_1$ such that $V(\xi(S)) > a$. Then there is some $T \in [T_1, S]$ such that $V(\xi(T)) = a$. We pick $T' \in [T_1, S]$ to be maximal with this property. It follows that $V(\xi(t)) > a$ for all $t \in (T', S]$. Applying the above argument with $\tau_1 = T'$ and $\tau_2 = S$, we have that $V(\xi(S)) \leq a = V(\xi(T'))$, a contradiction. So the claim holds.

We conclude that $V(\xi(t)) \leq \max\{a, V(\xi(0))\}$ for all $t \geq T_0$. Therefore the trajectory $\xi$ is precompact, and the first case in the Lemma is included in the second case, so we can assume that (9) holds. We claim that this means that $T_1 < \infty$, so that (8) holds with $T^* = T_1$. Indeed, if this were not true, then $\xi(t) \in C$ for all $t \geq T_0$, so $dV(\xi(t))/dt \leq -\alpha$ for almost all $t$, which gives $V(\xi(t)) \leq V(\xi(0)) - \alpha t$ for all $t \geq T_0$, which is impossible since $V$ is nonnegative.

To conclude that $\xi \to Z_V$, since $\xi$ is precompact we need only show that its omega-limit set $\Omega^+[\xi]$ is contained in $Z_V$. To see this, we pick any $z \in \Omega^+[\xi]$ and a sequence $\xi(t_i) \to z$. So $V(\xi(t_i)) \to V(z)$. If $z \notin Z_V$, let $a := V(z)/2 \neq 0$. Then Property (8) gives that $V(\xi(t_i)) \leq a$ for all $i$ large enough, which says that $\limsup_i V(\xi(t_i)) \leq a$, contradicting $V(\xi(t_i)) \to 2a$. Thus $z \in Z_V$. ∎

**Theorem 1** *Suppose given a behavior $\mathcal{R}$ of a system (6), and for for each compact subset $U \subseteq \mathcal{U}$, a $U$-decrease function $V_U$. Then:*

1. *For each compact subset $U \subseteq \mathcal{U}$ and for every $(\omega, \eta) \in \mathcal{R}$ satisfying $\omega \to U$, it holds that $\eta \to h(Z_{V_U})$.*

2. *If there is a compact subset $U_0 \subseteq \mathcal{U}$ and there is some $\gamma \in \mathcal{K}_\infty$ such that*

$$|h(Z_{V_U})| \leq \gamma(|U|) \tag{10}$$

   *for every compact $U \subseteq U_0$, then $\mathcal{R}$ has Cauchy gain $\gamma$ on $U_0$ and incremental limit gain $\gamma$ on $U_0$.*

3. *If there is some $\gamma \in \mathcal{K}_\infty$ such that (10) holds for every compact subset $U \subseteq \mathcal{U}$, then $\mathcal{R}$ has Cauchy gain $\gamma$ and incremental limit gain $\gamma$.*

*Proof.* We will first show that, for every triple $(\omega, \xi, \eta)$ with $\omega$ precompact and Lebesgue-measurable, such that $\dot{\xi}(t) = f(\xi(t), \omega(t))$ and $\eta(t) = h(\xi(t))$, the following properties hold:

(i) $\xi$ and $\eta$ are precompact;

(ii) for each compact $U \subseteq \mathcal{U}$, if $\omega \to U$ then $\xi \to Z_{V_U}$ and $\eta \to h(Z_{V_U})$.

(This will prove, in particular, the first assertion of the theorem.) Since $\omega$ is precompact, there exists some compact set, let us call it $U'$, such that $\omega(t) \in U'$ for all $t$. The first case in Lemma 4.2, applied to $U'$ and $V' = V_U$ with $T_0 = 0$, gives that $\xi(t) \to Z_{V'}$, so, being continuous as a function of $t$, $\xi$ is precompact. Next we apply once more Lemma 4.2, using now the second case with any given compact $U$ and the $U$-decrease function $V = V_U$, to conclude that $\xi \to Z_V$. Since the set $Z_V$ is compact and the mapping $h$ is continuous, it follows that also $\eta(t) = h(\xi(t)) \to h(Z_V)$ as $t \to \infty$, and $\eta$ is precompact as well.

Next, we pick any compact $U_0 \subseteq \mathcal{U}$ as in the second assertion, so that for every compact subset $U \subseteq U_0$ we have that $|h(Z_{V_U})| \leq \gamma(|U|)$, and let $\Gamma(U) := h(Z_{V_U})$. We will apply

Lemma 3.1. Property (c) in that Lemma holds by definition of $\Gamma$. Also, (a) holds, by (i). To prove that (b) is true, suppose that $\Omega^+[\omega] \subseteq U \subseteq U_0$. Since $\omega$ is precompact, $\omega \to \Omega^+[\omega]$, so also $\omega \to U$. By (ii), we know that $\eta \to h(Z_{V_U})$, and so $\Omega^+[\eta] \subseteq h(Z_{V_U}) = \Gamma(U)$, as desired. The Lemma then says that $\mathcal{R}$ has Cauchy gain $\gamma$ on $U_0$ and also has incremental limit gain $\gamma$ on $U_0$.

Finally, given an arbitrary pair $(\omega, \eta) \in \mathcal{R}$, we pick some compact $U_0$ such that $\omega(t) \in U_0$; then the just-shown Cauchy property on $U_0$ gives that $\|\eta\|_{\mathrm{aa}} \leq \gamma(\|\omega\|_{\mathrm{aa}})$, and similarly for the incremental gain. ∎

## 5 A Class of Examples

As an illustration of gain computations and small-gain stability arguments, we consider systems which consist of cascades of several subsystems, each of which can be individually described by some ordinary differential equation $\dot{x}_i = f(x_i, u_i)$ with input $u_i$. The input $u_1 = u$ to the first of the systems in the cascade is an external one, while the intermediate inputs $u_i$, $i > 1$, between two stages depend on the state of the preceding stage. In a biological application, $x_i(t)$ might represent the amount present, at any given time $t$, of the activated form $E_i^*$ of an enzyme $E_i$ whose production rate is, in turn, dependent upon the amount present of the activated form $E_{i-1}^*$ of the enzyme $E_{i-1}$. We allow transport delays in between stages. This leads to systems given by sets of delay-differential equations as follows:

$$\begin{aligned}
\dot{x}_1(t) &= f_1(x_1(t), u(t)) \\
\dot{x}_2(t) &= f_2(x_2(t), x_1(t - \tau_1)) \\
&\vdots \\
\dot{x}_n(t) &= f_n(x_n(t), x_{n-1}(t - \tau_{n-1}))
\end{aligned}$$

where $\tau_1, \ldots, \tau_{n-1} \geq 0$ are the delays among the stages of the process (the particular case in which there are no delays is included in this formalism by setting all $\tau_i = 0$). See Figure 3, where $x_i^\tau(t) := x_i(t - \tau)$ and we use $\mathcal{R}_i$ to denote the behavior associated to the system $\dot{x} = f_i(x, u)$ with output $y = x$. One often asks about such systems, see e.g. [3, 10, 17, 18], whether adding a feedback loop from the last stage to the first, as shown in Figure 4, might introduce instabilities, such as oscillations or even chaotic behavior. Specifically, one may have, for instance, that the action of $u$ on the first subsystem is inhibited by the final product $x_n$. Assuming that all the variables $x_i$ as well as the external input $u$ are scalar, and take only nonnegative values (as is for instance the case when variables represent chemical concentrations), a typical model for inhibition is obtained when the first equation becomes

$$\dot{x}_1(t) = f\left(x_1(t), \frac{u(t)}{1 + kx_n(t - \tau_n)}\right)$$

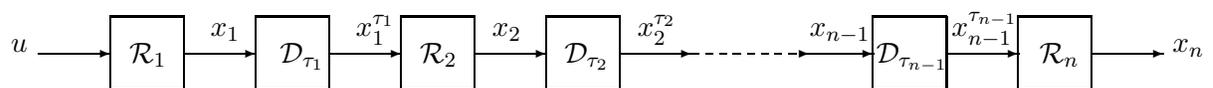

Figure 3: Cascade of $\mathcal{R}_i$'s and Delays

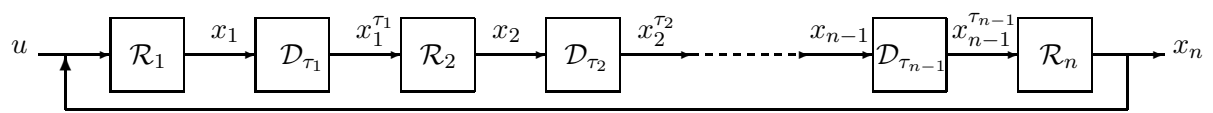

Figure 4: Cascade of $\mathcal{R}_i$'s and Delays, With Feedback to First Stage

(other expressions for inhibition are also possible, of course), where the "gain" $k \geq 0$ serves to parametrize the feedback strength. Note that we are also allowing for an additional delay in the feedback. Suppose that we are interested in analyzing the case in which $u(t)$ equals a constant value $\mu$, so that the effective input being fed to the first subsystem is $\omega(t) = \psi(x_n^{\tau_n}(t)) = \psi(x_n(t - \tau_n))$, where $\psi$ is the following function:

$$\psi(x) = \frac{\mu}{1 + kx} \qquad (11)$$

(for some particular values of $\mu$ and $k$). The closed-loop system may be viewed as the feedback interconnection $\mathcal{R} \cap \mathcal{M}_\psi^{-1}$ of the memoryless behavior $\mathcal{M}_\psi$ with the behavior of the forward composite system with output $x_n^{\tau_n}$, that is, the composition

$$\mathcal{R} = \mathcal{D}_{\tau_n} \circ \mathcal{R}_n \circ \mathcal{D}_{\tau_{n-1}} \circ \ldots \circ \mathcal{D}_{\tau_1} \circ \mathcal{R}_1.$$

Suppose that each $\mathcal{R}_i$ has Cauchy gain $\gamma_i$ and also incremental limit gain $\gamma_i$. Then $\mathcal{R}$ has gain $\gamma = \gamma_n \circ \ldots \circ \gamma_1$ of both types. Since $\psi$ has Lipschitz constant $k\mu$, $\mathcal{M}_\psi$ has both gains $k\mu I$. Therefore, provided that the small gain condition

$$\gamma(k\mu r) < r \quad \forall\, r > 0 \qquad (12)$$

holds, one concludes from Lemma 1.5 that there is some value $\bar{u}$ such that, for every solution of the closed-loop system, $\omega \to \bar{u}$. This means, in turn, because each $\mathcal{R}_i$ has a Cauchy gain, that every state variable $x_i$ converges to a unique equilibrium (independently of initial conditions).

The problem, therefore, is to estimate gains $\gamma_i$ for the systems $\mathcal{R}_i$. We briefly discuss one situation, itself of great interest, in which estimates can be obtained.

We suppose given intervals $\mathcal{X}_i = [a_i, b_i] \subseteq \mathbb{R}_{\geq 0}$, and sets $\mathcal{U}_i \subseteq \mathbb{R}_{\geq 0}$, $i = 1, \ldots, n$, with $\mathcal{U}_i = \mathbb{R}_{\geq 0}$ for $i > 1$, such that any solution of $\dot{x} = f_i(x, u)$ with initial conditions in $\mathcal{X}_i$ and input with values in $\mathcal{U}_i$ remains in $\mathcal{X}_i$. In a typical biological application, one may have $\mathcal{X}_i = [0, x_i^{\max}]$, where $x_i^{\max}$ is the maximum possible amount of a substance, such as the activated form of an enzyme, that may be synthesized. Also, suppose given, for each $i$, a strictly increasing and onto function

$$g_i : [a_i, b_i) \to \mathbb{R}_{\geq 0}$$

with the following properties:

1. the restriction of $g_i^{-1} : \mathbb{R}_{\geq 0} \to [a_i, b_i)$ to $\mathcal{U}_i$ is Lipschitz, with Lipschitz constant $\lambda_i$;

2. $x < g_i^{-1}(u) \Rightarrow f_i(x, u) > 0$, and $x > g_i^{-1}(u) \Rightarrow f_i(x, u) < 0$, for every $u \in \mathcal{U}_i$ and $x \in \mathcal{X}_i$.

Then, $\mathcal{R}$ admits both Cauchy and incremental limit gain $\gamma(r) = \lambda_1 \ldots \lambda_n r$, and (12) becomes:

$$k\mu < \frac{1}{\lambda_1 \ldots \lambda_n}. \qquad (13)$$

This is proved as follows. Fixing any $i = 1, \ldots, n$, we drop subscripts and write $\mathcal{X} = \mathcal{X}_i, \mathcal{U} = \mathcal{U}_i$, $g = g_i, \lambda = \lambda_i$. For each compact subset $U \subseteq \mathcal{U}$, we let $c := \min U$ and $d := \max U$, so that $U \subseteq [c, d]$ and $|U| = d - c$, and define $V_U : \mathcal{X} \to \mathbb{R}_{\geq 0}$ to be the distance from any $x$ to the set

$$g^{-1}([c,d]) \;=\; [g^{-1}(c), g^{-1}(d)]$$

(recall that $g^{-1}$ is an increasing function). So $Z_{V_U} = [g^{-1}(c), g^{-1}(d)]$, and $V_U(x) = g^{-1}(c) - x$ if $g(x) < c$, $V_U(x) = x - g^{-1}(d)$ if $g(x) > d$. Note that $Z_{V_U}$ has diameter $|Z_{V_U}| \leq \lambda |d-c| = \lambda |U|$. Furthermore, $V_U$ is a $U$-decrease function, because $V_U$ is continuous, proper, differentiable outside $Z_{V_U}$, and satisfies the decrease condition. Indeed, pick any $x \in \mathcal{X} \setminus Z_{V_U}$ and $u \in U$. There are two cases: $x < g^{-1}(c)$ or $x > g^{-1}(d)$. In the first case, $x < g^{-1}(c) \leq g^{-1}(u)$ (since $u \in [c,d]$), so $\nabla V(x) \cdot f(x,u) = (-1) \cdot f(x,u) < 0$, and the second case is similar. Theorem 1 (part 3) then applies, so the corresponding behavior $\mathcal{R}_i$ has both Cauchy and incremental limit gain $\gamma_i(r) = \lambda_i r$, and therefore the cascade (with arbitrary delays) has both gains equal to $\gamma$, as claimed.

The requirement that each $g_i^{-1} : \mathcal{U}_i \to [a_i, b_i)$ must have Lipschitz constant $\lambda_i$ can be relaxed to:

$$\text{the restriction of } g_i^{-1} \text{ to } U_i \text{ is Lipschitz with Lipschitz constant } \lambda_i \tag{14}$$

for each $i = 1, \ldots, n$, where we define, inductively, the intervals $U_1 := \mathcal{U}_1$ and

$$U_{i+1} \;:=\; Z_{V_{U_i}} \;=\; g_i^{-1}(U_i) \;\subseteq\; \mathcal{U}_{i+1}$$

for $i = 1, \ldots, n$. For each compact subset $U \subseteq U_i$, we have that $|Z_{V_U}| \leq \lambda_i |U|$, so, by Theorem 1 (part 2, applied $n$ times, with $U_0 = U_i$, $i = 1, \ldots, n$) we conclude that each $\mathcal{R}_i$ has Cauchy gain as well as incremental limit gain $\gamma_i(r) = \lambda r$ on $U_i$. In addition, for each $i = 1, \ldots, n-1$ and every $(\omega, \eta) \in \mathcal{R}_i$ such that $\omega \to U_i$, Theorem 1 (part 1) insures that $\eta \to U_{i+1}$. The pure delays $D_{\tau_i}$ have identity gains, and of course satisfy $\omega \to U_{i+1} \Rightarrow \eta \to U_{i+1}$. So, arguing (for each consecutive pair of subsystems in the cascade) as in Section 2.1, we conclude that $\mathcal{R}$ *admits both Cauchy and incremental limit gain $\gamma(r) = \lambda_1 \ldots \lambda_n r$ provided condition (14) holds.*

Let us specialize even further. We assume from now on that each function $f_i$ has the following form:

$$f_i(x, u) \;=\; -\alpha_i(x) + u\beta_i(x),$$

where the functions $\alpha_i$ and $\beta_i$ are nonnegative on $\mathcal{X} = [a_i, b_i]$, $\alpha_i$ is strictly increasing and $\beta_i$ is strictly decreasing, and $\alpha_i(a_i) = \beta_i(b_i) = 0$. Since $f_i(a_i, u) > 0$ and $f_i(b_i, u) < 0$ for all $u \in \mathcal{U}_i$, the interval $\mathcal{X}_i = [a_i, b_i]$ is invariant. For each $i$, we let

$$g_i(x) \;:=\; \frac{\alpha_i(x)}{\beta_i(x)}$$

for $x \in [a_i, b_i)$. Then $g_i$ is strictly increasing, because $\alpha_i$ increases and $\beta_i$ decreases, and it satisfies $g_i(a_i) = 0$ and $g_i(x) \to \infty$ as $x \to b^-$, so it is onto $\mathbb{R}_{\geq 0}$. Given any $x \in [a_i, b_i)$ and $u \in \mathcal{U}_i$, if $x < g_i^{-1}(u)$ then $\alpha_i(x)/\beta_i(x) = g_i(x) < u$ implies $f_i(x, u) > 0$ and similarly $x > g_i^{-1}(u) \Rightarrow f_i(x, u) < 0$. If $x = b_i$ then $f_i(x, u) = -\alpha_i(b) < 0$. In conclusion, the functions $g_i(x) := \frac{\alpha_i(x)}{\beta_i(x)}$ are as required for the above computations, and the gains can be computed in terms of the Lipschitz constants in (14), computed on the sets $U_i$.

## 5.1 MAPK Cascades

As an application, and in fact the original motivation for this study, we pick the case when every $\mathcal{X}_i = [0,1]$ (for now, we let $\mathcal{U}_1$ be arbitrary, but it will be restricted below), and

$$\alpha_i(x) \;=\; \frac{b_i x}{c_i + x} \quad \text{and} \quad \beta_i(x) \;=\; \frac{d_i(1-x)}{e_i + (1-x)}$$

(for some positive constants $b_i, c_i, d_i, e_i$). We wish to study the stability of the inhibitory closed-loop system:

$$\dot{x}_1 \;=\; -\frac{b_1 x_1}{c_1 + x_1} + \frac{\mu}{1 + k x_3^{\tau_3}} \frac{d_1(1-x_1)}{e_1 + (1-x_1)} \tag{15}$$

$$\dot{x}_2 \;=\; -\frac{b_2 x_2}{c_2 + x_2} + x_1^{\tau_1} \frac{d_2(1-x_2)}{e_2 + (1-x_2)} \tag{16}$$

$$\dot{x}_3 \;=\; -\frac{b_3 x_3}{c_3 + x_3} + x_2^{\tau_2} \frac{d_3(1-x_3)}{e_3 + (1-x_3)} \tag{17}$$

where the $\tau_i \geq 0$. Such equations arise, for instance, as follows.

*Mitogen-activated protein kinase (MAPK) cascades* constitute a type of "biological module" or "subcircuit" which is implicated, in several variants, in a large variety of eukaryotic cell signal transduction processes, cf. [4, 6, 19]. This highly conserved signaling cascade processes inputs – themselves triggered, in turn, by extracellular stimuli – into output signals responsible for diverse cellular behaviors: proliferation, growth, and differentiation (hence the name), as well as movement, stress responses, and death. The basic mechanism is that of a cascade of three subsystems, each of which consists of one or more reversible enzyme activations.

Subject to conservation laws (stoichiometry relations), each subsystem is usually described in terms of one or two differential equations, see for instance [4, 8, 13, 14]. (The models in these references do not include delays among levels in the cascade, but in our formalism, arbitrary delays do not affect the results, and they are biologically plausible.) For simplicity, we pick a model with one equation in each level, as in [13, 14] and the last section of [8]. In general, one would have all $\mathcal{X}_i = [0, E_i]$, where $E_i$ is the total amount of enzyme present (activated plus nonactivated), but we may take all $\mathcal{X}_i = [0, 1]$ after nondimensionalization, as was done in [13, 14]. These references, and especially Kholodenko's work, emphasized the fact that inhibitory feedback, which seems to be present in naturally occuring cells, could, theoretically, produce oscillations. As oscillations in MAPK cascades do not appear to occur naturally, an interesting mathematical question is to find conditions on gains insuring lack of oscillations. We do this next.

We first show that the functions $g_i^{-1}$ are always Lipschitz, and compute an estimate of the constant. This will prove that stability is preserved under small enough feedback, but the estimates may be too conservative. After that, we describe a numerical approach which allows sometimes quite tight estimates. To obtain the general bounds, we compute, for $x \in [0, 1)$:

$$g_i'(x) \;=\; \frac{b_i}{d_i} \frac{e_i x^2 + c_i(x-1)^2 + c_i e_i}{(c_i + x)^2 (1-x)^2}$$

and use (taking derivatives to minimize) these lower bounds:

$$e_i x^2 + c_i(x-1)^2 \geq \frac{c_i e_i}{e_i + c_i}$$

$$(c_i + x)(1 - x) \leq (1/4)(c_i + 1)^2$$

so as to obtain the following estimate:

$$g_i'(x) \geq \delta_i := 16 \frac{b_i}{d_i} \frac{c_i e_i}{(c_i + 1)^4} \left(1 + \frac{1}{e_i + c_i}\right) > 0.$$

Thus $\lambda_i = 1/\delta_i$ are Lipschitz constants as desired.

Since we are ultimately interested in the effect of inhibitory feedback given by a function $\psi(\xi) = \frac{\mu}{1+k\xi}$ as in (11), where $\xi$ is a (possibly delayed) value of $x_3$, which ranges over the interval $[0, 1]$, the only inputs $u$ that must be considered are those in $[\frac{\mu}{1+k}, \mu]$, So, from now, on we suppose that the first input set has the form:

$$\mathcal{U}_1 = [\bar{u}, \infty)$$

for some $\bar{u} > 0$. The appropriate value of $\bar{u}$ will depend upon prior upper bounds on the feedback gains $k$ and external inputs $\mu$ which are of interest. We let $U_1 = \mathcal{U}_1$ and introduce $U_{i+1} = g_i^{-1}(U_i)$ for $i = 1, 2, 3$ as done earlier, and compute a Lipschitz constant $\lambda_i$ for the restriction of $g_i^{-1}$ to $U_i$ for each $i = 1, 2, 3$ as in (14). Writing $\bar{x}_0 = \bar{u}$ and $\bar{g}_i(x_i) := \bar{x}_{i-1}$ for $i = 1, 2, 3$, we have that $U_2 \subseteq [\bar{x}_1, 1]$, $U_3 \subseteq [\bar{x}_2, 1]$, $U_4 \subseteq [\bar{x}_3, 1]$. When applying the small-gain condition (13), the quantity of interest is the "total gain" $\lambda = \lambda_1 \lambda_2 \lambda_3$ of the cascade. An upper bound on each $\lambda_i$ can be obtained, for any given $\bar{u}$, by maximizing the derivative of $g_i^{-1}$ on $U_i$, or equivalently (and most conveniently, as this involves no functional inversion and each $g_i$ is a simple rational function) as $\lambda_i = \theta_i(\bar{u}_i)^{-1}$, where $\theta_i(\bar{u})$ is the *minimum of $g_i'$ on $[\bar{x}_i, 1]$*. The functions $\theta_i$ depend, of course, on the actual parameters $b_i$, etc. Defining $\theta(\bar{u}) := \theta_1(\bar{u})\theta_2(\bar{u})\theta_3(\bar{u})$, we conclude that *every solution of the delay-differential system (15)-(16)-(17) will satisfy $x_i(t) \to 0$ as $t \to \infty$ provided that the following condition holds*:

$$k < \min\left\{\frac{\theta(\bar{u})}{\mu}, \frac{\mu}{\bar{u}} - 1\right\}. \tag{18}$$

The first term represents the small-gain condition (13) (since $\bar{u} = \frac{1}{\lambda}$), while the second insures that $u = \psi(x_3^{\tau_3}) = \mu(1 + kx_3^{\tau_3})^{-1}$ belongs to the input set $\mathcal{U}_1 = [\bar{u}, \infty)$ for all $x_3 \in [0, 1]$. Figure 5 provides a plot (obtained using the "minimize" function in Maple) of the function $\theta$

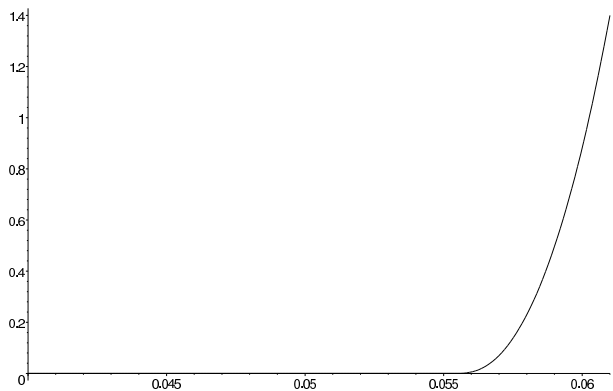

Figure 5: The Function $\theta$

for a concrete numerical example which takes the coefficients used in [13]:

$$b_1 = c_1 = e_1 = b_2 = 0.1, c_2 = e_2 = c_3 = e_3 = 0.01, b_3 = 0.5, d_1 = d_2 = d_3 = 1.\qquad(19)$$

There is a sharp transition at around $\bar{u} = 0.055$. In order to obtain $\lambda$ small, let us say $\lambda < 1$, it is reasonable to pick the value $\bar{u} = 0.061$, for which $\theta(\bar{u}) \approx 1.39933$ and so $\lambda = 1/\theta(\bar{u}) \approx 0.71463$ (picking instead $\bar{u} = 0.06$ would give $\lambda \approx 1.134$). If we now pick $\mu = 0.3$, as in [13], the constraints (18) become $k < \min\{4.6644, 3.918\} = 3.918$. In conclusion, *any feedback with gain $k < 3.9$ is guaranteed to preserve stability* (for external input $u(t) \equiv 0.3$, and for any delays).

One may ask how tight is the bound $k < 3.9$. For the system with no delays, numerical experimentation leads to the conclusion that, as the parameter $k$ increases, a Hopf bifurcation occurs at around $k = 5.1$. The solution with initial conditions $x(0) = 0$, and $\mu = 0.3$, $k = 5.2$ is oscillatory and is plotted in Figure 6. (Linear instability in this context always leads to existence of periodic solutions, cf. [3, 18].) Thus the bound is fairly tight even when there are no delays.

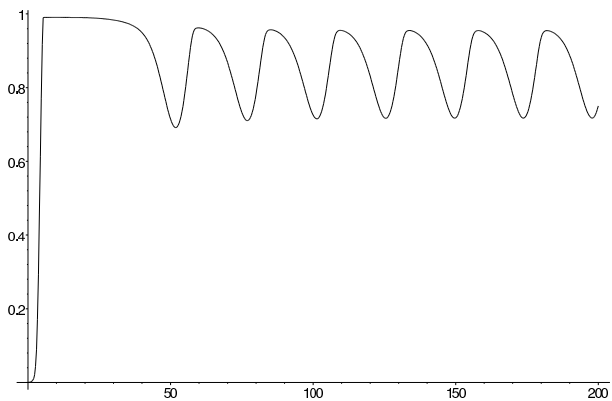

Figure 6: Oscillations when $k = 5.2$ and $\mu = 0.3$

## 5.2 Local Exponential Stability

For systems of the special form being considered, linearized stability analysis can be employed in order to determine local exponential stability of the feedback system, appealing to classical small-gain theorems as in [1, 11, 12, 21]. (A computational difficulty when doing this, however, is that feedback laws such as $u = \psi(x)$ change equilibria, and the new equilibria are usually not computable in closed-form.) Given any system of the form $\dot{x} = f(x,u) = -\alpha(x) + u\beta(x)$ as above, and any $\bar{u}$, we may compute the linearized system $\dot{z} = az + bv$ at the equilibrium values $x = \bar{x} = g^{-1}(\bar{u})$ and $u = \bar{u}$. This has $a = \frac{\partial f}{\partial x}(\bar{x}, \bar{u})$ and $b = \frac{\partial f}{\partial u}(\bar{x}, \bar{u})$. From $a/b = (-\alpha'(\bar{x}) + \bar{u}\beta'(\bar{x}))/\beta(\bar{x})$ and $\bar{u} = \alpha(\bar{x})/\beta(\bar{x}) = g(\bar{x})$, one obtains:

$$\left|\frac{a}{b}\right| = \frac{\alpha'(\bar{x})\beta(\bar{x}) - \alpha(\bar{x})\beta'(\bar{x})}{\beta^2(\bar{x})} = g'(\bar{x}) = \frac{1}{(g^{-1})'(\bar{u})}.$$

On the other hand, the $H_\infty$ gain (induced $\mathcal{L}^2$ operator norm) of the system $\dot{z} = az + bv$ is $\max\{b/|i\omega + a|, \omega \in \mathbb{R}\} = |b/a|$. We conclude that $(g^{-1})'(\bar{u})$ equals the $H_\infty$ gain of the linearized system. By induction on the cascade, *the linearized cascade system has an $H_\infty$ gain equal to $\lambda$*, so the same inequalities on $\kappa$ and $\mu$ insure local exponential stability. (Linearized analysis cannot by itself guarantee the existence of equilibria for the closed-loop system; only after equilibria are known to exist, this linearized small gain argument can be used.)

The "secant condition" for stability, see [2, 8, 10, 17, 18], is often used for linearized stability analysis of inhibitory feedback systems such as the ones considered here. This condition says that a matrix

$$\begin{pmatrix} \alpha_1 & 0 & \cdots & 0 & -\beta_1 \\ \beta_2 & \alpha_2 & \cdots & 0 & -0 \\ \vdots & \vdots & & & \vdots \\ 0 & 0 & \cdots & \beta_n & \alpha_n \end{pmatrix}$$

with all $\alpha_i < 0$ and all $\beta_i > 0$ is Hurwitz provided that:

$$\left| \frac{\beta_1 \ldots \beta_n}{\alpha_1 \ldots \alpha_n} \right| < \left( \sec \frac{\pi}{n} \right)^n .$$

For $n = 3$ as in our example, this means that a closed-loop gain margin of $(\sec \pi/3)^3 = 8$ can be tolerated while preserving stability, for the linearized system (no assertion is made about phase, so no delays are allowed). The condition $k \leq \mu/\bar{u} - 1$, which guarantees that $u = \psi(x_3)$ belongs to the input set $[\bar{u}, \infty)$, must still be satisfied, but the bound (18) may be relaxed to $k < \min\{\frac{8\theta(\bar{u})}{\mu}, \frac{\mu}{\bar{u}} - 1\}$. In order to exploit this condition for the above numerical example, a bit of experimentation leads us to pick $\bar{u} = 0.05763$ (instead of 0.061), so $\lambda = 1/\theta(\bar{u}) \approx 6.32$. Using again $\mu = 0.3$, the condition $k \leq \mu/\bar{u} - 1$ imposes a constraint of approximately $k \leq 4.2$. For any such $k$, we have that $k \cdot \mu \cdot \lambda < 7.9632 < 8$, as wanted. In other words, the secant rule allows us to conclude that a feedback gain of $k \leq 4.2$ may be tolerated so as to guarantee linearized stability. This is only a slight improvement over the estimate $k \leq 3.9$ obtained from the nonlinear small-gain result, and it comes at the cost of assuming *no delays* and insuring only *local* stability for the original system.

# References

[1] Desoer, C.A., and M. Vidyasagar, *Feedback Synthesis: Input-Output Properties*, Academic Press, New York, 1975.

[2] Dibrov, B.F., A.M. Zhabotinsky. and B.N. Kholodenko, "Dynamic stability of steady states and static stabilization in unbranched metabolic pathways," *J. Math. Biol.* **15**(1982): 51–63.

[3] Hastings, S., J. Tyson, and D. Webster, "Existence of periodic solutions for negative feedback cellular control systems," *J. Diff. Eq.* **25**(1977): 39–64.

[4] Huang, C.-Y.F., and J.E. Ferrell, "Ultrasensitivity in the mitogen-activated protein kinase cascade," *Proc. Natl. Acad. Sci. USA* **93**(1996): 10078–10083.

[5] Ingalls, B., and E.D. Sontag, "A small-gain lemma with applications to input/output systems, incremental stability, detectability, and interconnections," submitted, 2001.

[6] Lauffenburger, D.A., "A computational study of feedback effects on signal dynamics in a mitogen-activated protein kinase (MAPK) pathway model" *Biotechnol. Prog.* **17**(2001): 227–239.

[7] Jiang, Z.-P., A. Teel, and L. Praly, "Small-gain theorem for ISS systems and applications," *Mathematics of Control, Signals, and Systems* **7**(1994): 95-120.

[8] Kholodenko, B.N., "Negative feedback and ultrasensitivity can bring about oscillations in the mitogen-activated protein kinase cascades," *Eur. J. Biochem* **267**(2000): 1583–1588.

[9] Mareels, I.M.Y., and D.J. Hill, "Monotone stability of nonlinear feedback systems," *J. Math. Systems, Estimation and Control* **2**(1992): 275–291.

[10] Othmer, H.G., "The qualitative dynamics of a class of biochemical control circuits," *J. Math. Biol.* **3**(1976): 53–78.

[11] Safonov, M., *Stability and Robustness of Multivariable Feedback Systems*, The MIT Press, Cambridge, Massachusetts, 1980.

[12] Sandberg, I.W., "On the $L_2$-boundedness of solutions of nonlinear functional equations," *Bell System Technical Journal* **43**(1964): 1581–1599.

[13] Shvartsman, S.Y., H.S. Wiley, and D.A. Lauffenburger, "Autocrine loop as a module for bidirectional and context-dependent cell signaling," preprint, MIT Chemical Engineering department, December 2000.

[14] Shvartsman, S.Y., M.P. Hagan, A. Yacoub, P. Dent, H.S. Wiley, and D.A. Lauffenburger, "Context-dependent signaling in autocrine loops with positive feedback: Modeling and experiments in the EGFR system," *Am. J. Physiol. Cell Physiol.*, to appear.

[15] Sontag, E.D., "Smooth stabilization implies coprime factorization," *IEEE Trans. Automatic Control* **34**(1989): 435-443.

[16] Sontag, E.D., *Mathematical Control Theory, Deterministic Finite Dimensional Systems* (Second Edition), Springer-Verlag, New York, 1988.

[17] Thron, C.D., "The secant condition for instability in biochemical feedback-control .1. the role of cooperativity and saturability," *Bull. Math. Biology* **53**(1991): 383–401.

[18] Tyson, J.J., and H.G. Othmer, "The dynamics of feedback control circuits in biochemical pathways," in *Progress in Theoretical Biology (R. Rosen & F.M. Snell, eds.)*, Vol. 5, Academic Press, New York, 1978, pp. 1–62.

[19] Widmann, C., G. Spencer, M.B. Jarpe, and G.L. Johnson, "Mitogen-activated protein kinase: Conservation of a three-kinase module from yeast to human," *Physiol. Rev.* **79**(1999): 143–180.

[20] Willems, J.C., "On interconnections, control, and feedback," *IEEE Trans. Autom. Control* **42**(1997): 326–339.

[21] Zames, G., "On the input-output stability of time-varying nonlinear feedback systems. Part I: Conditions using concepts of loop gain, conicity, and positivity," *IEEE Trans. Autom. Control* **11**(1966): 228–238.